\documentclass[12pt,reqno]{amsart}

\usepackage{amssymb}

\newcommand{\kap}{\kappa}
\newcommand{\lam}{\lambda}
\renewcommand{\phi}{\varphi}

\newcommand{\Gam}{\Gamma}

\newcommand{\Z}{{\mathbb Z}}

\newcommand{\cH}{{\mathcal H}}
\newcommand{\cL}{{\mathcal L}}

\DeclareMathOperator{\Cay}{Cay}

\newcommand{\CP}[1][G]{\Cay_{#1}^+}

\newcommand{\pH}[1][H]{\phi_{#1}}

\newtheorem{lemma}{Lemma}
\newtheorem{theorem}{Theorem}
\newtheorem{corollary}{Corollary}
\newtheorem{proposition}{Proposition}

\newtheorem{theirtheorem}{Theorem}
\newtheorem{theirproposition}{Proposition}

\newcommand{\seq}{\subseteq}
\newcommand{\stm}{\setminus}
\newcommand{\est}{\varnothing}

\newcommand{\<}{\langle}
\renewcommand{\>}{\rangle}

\newcommand{\longc}{,\dotsc,}

\newcommand{\reft}[1]{~\ref{t:#1}}
\newcommand{\refl}[1]{~\ref{l:#1}}
\newcommand{\refc}[1]{~\ref{c:#1}}
\newcommand{\refp}[1]{~\ref{p:#1}}
\newcommand{\refs}[1]{~\ref{s:#1}}
\newcommand{\refb}[1]{~\cite{b:#1}}
\newcommand{\refe}[1]{~\eqref{e:#1}}

\newcommand{\be}{\begin{equation}}
\newcommand{\ee}{\end{equation}}
\author{David Grynkiewicz}
\address{Departament de Matem\`atica Aplicada IV,
  Universitat Polit\`ecnica de Catalunya,
  Campus Nord Edifici C3,
  Jordi Girona Salgado 1-3,
  Barcelona, Catalonia E-08034, Spain.}
\thanks{The first author is supported in part by the National Science
  Foundation, as an MPS-DRF postdoctoral fellow, under grant DMS-0502193.}
 \email{diambri@hotmail.com}

\author{Vsevolod F. Lev}
\address{Department of Mathematics, University of Haifa at Oranim,
  Tivon 36006, Israel}
\email{seva@math.haifa.ac.il}

\author{Oriol Serra}
\address{Departament de Matem\`atica Aplicada IV,
  Universitat Polit\`ecnica de Catalunya,
  Campus Nord Edifici C3,
  Jordi Girona Salgado 1-3,
  Barcelona, Catalonia E-08034, Spain}
\email{oserra@ma4.upc.edu}

\title{Connectivity of addition Cayley graphs}

\begin{document}
\baselineskip 16pt

\begin{abstract}
For any finite abelian group $G$ and any subset $S\seq G$, we determine the
connectivity of the addition Cayley graph induced by $S$ on $G$. Moreover, we
show that if this graph is not complete, then it possesses a minimum vertex
cut of a special, explicitly described form.
\end{abstract}

\maketitle

\section{Background: addition Cayley graphs}\label{s:background}

For a subset $S$ of the abelian group $G$, we denote by $\CP(S)$ the
addition Cayley graph induced by $S$ on $G$; recall that this is the
undirected graph with the vertex set $G$ and the edge set
$\{(g_1,g_2)\in G\times G\colon g_1+g_2\in S\}$. Note that $S$ is
not assumed to be symmetric, and that if $S$ is finite, then
$\CP(S)$ is regular of degree $|S|$ (if one considers each loop to
contribute $1$ to the degree of the corresponding vertex).

The twins of the usual Cayley graphs, addition Cayley graphs (also called
\emph{sum graphs}) received much less attention in the literature; indeed,
\refb{a} (independence number), \refb{cgw} and \refb{l2} (hamiltonicity),
\refb{chu} (expander properties), and \refb{g} (clique number) is a nearly
complete list of papers, known to us, where addition Cayley graphs are
addressed. To some extent, this situation may be explained by the fact that
addition Cayley graphs are rather difficult to study. For instance, it is
well-known and easy to prove that any connected Cayley graph on a finite
abelian group with at least three elements is hamiltonian, see \refb{mr};
however, apart from the results of \refb{cgw}, nothing seems to be known on
hamiltonicity of \emph{addition} Cayley graphs on finite abelian groups.
Similarly, the connectivity of a Cayley graph on a finite abelian group is
easy to determine, while determining the connectivity of an \emph{addition}
Cayley graph is a non-trivial problem, to the solution of which the present
paper is devoted. The reader will see that investigating this problem leads
to studying rather involved combinatorial properties of the underlying group.

\section{Preliminaries and summary of results}\label{s:intro}

Let $\Gam$ be a graph on the finite set $V$. The (vertex) connectivity of
$\Gam$, denoted by $\kap(\Gam)$, is the smallest number of vertices which are
to be removed from $V$ so that the resulting graph is either disconnected or
has only one vertex. Clearly, if $\Gam$ is complete, then $\kap(\Gam)=|V|-1$,
while otherwise we have $\kap(\Gam)\le|V|-2$, and $\kap(\Gam)$ can be
alternatively defined as the size of a minimum vertex cut of $\Gam$. (A
complete graph does not have vertex cuts.) Evidently, vertex cuts and
connectivity of a graph are not affected by adding or removing loops.

Our goal is to determine the connectivity of the addition Cayley graphs,
induced on finite abelian groups by their subsets, and accordingly we use
additive notation for the group operation. In particular, for subsets $A$
and $B$ of an abelian group, we write
  $$ A\pm B:=\{a\pm b\colon a\in A,\,b\in B\}, $$
which is abbreviated by $A\pm b$ in the case where $B=\{b\}$ is a singleton
subset.

For the rest of this section, we assume that $S$ is a subset of the finite
abelian group $G$.

It is immediate from the definition that, for a subset $A\seq G$, the
neighborhood of $A$ in $\CP(S)$ is the set $S-A$, and it is easy to derive
that $\CP(S)$ is complete if and only if either $S=G$, or $S=G\stm\{0\}$ and
$G$ is an elementary abelian $2$-group (possibly of zero rank). Furthermore,
it is not difficult to see that $\CP(S)$ is connected if and only if $S$ is
not contained in a coset of a proper subgroup of $G$, with the possible
exception of the non-zero coset of a subgroup of index $2$; this is
\cite[Proposition~1]{b:l2}. Also, since $\CP(S)$ is $|S|$-regular, we have
the trivial bound $\kap(\CP(S))\le|S|$.

If $H$ is a subgroup of $G$ and $g$ is an element of $G$ with $2g\in S+H$,
then $g+H\seq S-(g+H)$; consequently, the boundary of $g+H$ in $\CP(S)$ has
size
  $$ |(S-(g+H))\stm(g+H)| = |S+H|-|H|. $$
Assuming in addition that $S+H\neq G$, we obtain $(S-(g+H))\cup(g+H)=S+H-g\ne
G$, implying $\kap(\CP(S))\le|S+H|-|H|$. Set
  $$ 2\ast G := \{ 2g\colon g\in G \}, $$
so that the existence of $g\in G$ with $2g\in S+H$ is equivalent to the
condition $(S+2\ast G)\cap H\neq\est$. Motivated by the above observation, we
define
  $$ \cH_G(S) := \{ H\le G\colon (S+2\ast G)\cap H\neq\est,\, S+H\neq G \} $$
and let
  $$ \eta_G(S):= \min \{ |S+H|-|H| \colon H\in\cH_G(S) \}. $$
In the latter definition and throughout, we assume that the minimum of an
empty set is infinite, and we allow comparison between infinity and real
numbers according to the ``na{\i}ve'' rule. Thus, for instance, we have
$\kap(\CP(S))\le\eta_G(S)$ even if $\cH_G(S)$ is vacuous.

Another important family of sets with small boundary is obtained as
follows. Suppose that the subgroups $L\le G_0\le G$ and the element
$g_0\in G_0$ satisfy
\begin{itemize}
\item[(i)]  $|G_0/L|$ is even and larger than $2$;
\item[(ii)] $S+L=(G\stm G_0)\cup(g_0+L)$.
\end{itemize}
Fix $g\in G_0\stm L$ with $2g\in L$ and consider the set
$A:=(g+L)\cup(g+g_0+L)$. The neighborhood of this set in $\CP(S)$ is
  $$ S-A = (G\stm G_0) \cup (g+L) \cup (g+g_0+L) = (G\stm G_0) \cup A, $$
whence  $(S-A)\cup A\neq G$ and $|(S-A)\stm A|=|G\stm
G_0|=|S+L|-|L|$. Consequently, $\kap(\CP(S))\le|S+L|-|L|$. With this
construction in mind, we define $\cL_G(S)$ to be the family of all
those subgroups $L\le G$ for which a subgroup $G_0\le G$, lying
above $L$, and an element $g_0\in G_0$ can be found so that
properties (i) and (ii) hold, and we let
  $$ \lam_G(S):= \min \{ |S+L|-|L| \colon L\in\cL_G(S) \}. $$
Thus, $\kap(\CP(S))\le\lam_G(S)$.

Our first principal result is the following.
\begin{theorem}\label{t:main1}
If $S$ is a proper subset of the finite abelian group $G$, then
  $$ \kap(\CP(S)) = \min \{ \eta_G(S), \lam_G(S), |S| \}. $$
\end{theorem}

Let $\Gam$ be a graph on the vertex set $V$. We say that the
non-empty subset $V_0\subset V$ is a \emph{fragment} of $\Gam$ if
the neighborhood $N(V_0)$ of $V_0$ satisfies $|N(V_0)\stm
V_0|=\kap(\Gam)$ and $N(V_0)\cup V_0\ne V$; that is, the boundary of
$V_0$ is a minimum vertex-cut, separating $V_0$ from the (non-empty)
remainder of the graph. Notice that if $\Gam$ is not complete, then
it has fragments; for instance, if $\Gam'$ is obtained from $\Gam$
by removing a minimum vertex cut, then the set of vertices of any
connected component of $\Gam'$ is a fragment of $\Gam$.

As the discussion above shows, if $\kap(\CP(S))=\eta_G(S)$, then $\CP(S)$ has
a fragment which is a coset of a subgroup $H\in\cH_G(S)$ with
$|S+H|-|H|=\eta_G(S)$; similarly, if $\kap(\CP(S))=\lam_G(S)$, then $\CP(S)$
has a fragment which is a union of at most two cosets of a subgroup
$L\in\cL_G(S)$ with $|S+L|-|L|=\lam_G(S)$.

The reader will easily verify that Theorem \reft{main1} is an
immediate corollary of Theorem \reft{main2} below. The latter shows
that the minimum in the statement of Theorem \reft{main1} is
attained, with just one exception, on either $\eta_G(S)$ or $|S|$.
Being much subtler, Theorem \reft{main2} is also more technical, and
to state it we have to bring into consideration a special sub-family
of $\cL_G(S)$. Specifically, let $\cL^\ast_G(S)$ be the family of
those subgroups $L\le G$ such that for some $G_0\le G$, lying above
$L$, and some $g_0\in G_0$, the following conditions hold:
\begin{itemize}
\item[(L1)] $G_0/L$ is a cyclic $2$-group of order $|G_0/L|\ge 4$, and
  $\<g_0\>+L=G_0$;
\item[(L2)] $G/G_0$ is an elementary abelian $2$-group (possibly of zero
  rank);
\item[(L3)] $\exp(G/L)=\exp(G_0/L)$;
\item[(L4)] $S+L=(G\stm G_0)\cup(g_0+L)$ and $S\cap(g_0+L)$ is not
  contained in a proper coset of $L$.
\end{itemize}

A little meditation shows that $\cL^\ast_G(S)\seq\cL_G(S)$ and that
conditions (L1)--(L3) imply
  $$ G/L \cong (G_0/L)\oplus(\Z/2\Z)^r \cong (\Z/2^k\Z)\oplus(\Z/2\Z)^r, $$
for some $k\ge 2$ and $r\ge 0$. Notice also that if $L$, $G_0$, and $g_0$ are
as in (L1)--(L4), and $G_0=G$, then $L$ is a subgroup of $G$ of index at
least $4$, and $S$ is contained in an $L$-coset, whence $\CP(S)$ is
disconnected.

\begin{theorem}\label{t:main2}
Let $S$ be a proper subset of the finite abelian group $G$. There exists
at most one subgroup $L\in\cL^\ast_G(S)$ with $|S+L|-|L|\le |S|-1$.
Moreover,
\begin{itemize}
\item[(i)] if $L$ is such a subgroup, then
  $\kap(\CP(S))=\lam_G(S)=|S+L|-|L|$ and
  $\eta_G(S)\ge |S|$;
\item[(ii)] if such a subgroup does not exist, then
  $\kap(\CP(S))=\min\{\eta_G(S),|S|\}$.
\end{itemize}
\end{theorem}

Postponing the proof to Section \refs{proofs}, we now list some of the
consequences.

\begin{corollary}
Let $S$ be a proper subset of the finite abelian group $G$ such that $\CP(S)$
is connected. If either $|S|\le|G|/2$ or $G$ does not contain a subgroup
isomorphic to $(\Z/4\Z)\oplus(\Z/2\Z)$, then
$\kap(\CP(S))=\min\{\eta_G(S),|S|\}$.
\end{corollary}

\begin{proof}
If $\kap(\CP(S))\ne\min\{\eta_G(S),|S|\}$, then by Theorem \reft{main2}
there exists $L\in\cL^\ast_G(S)$ with $|S+L|-|L|\le|S|-1$. Choose $L\le
G_0\le G$ and $g_0\in G_0$ satisfying (L1)--(L4). Since $\CP(S)$ is
connected, the subgroup $G_0$ is proper. Consequently,
  $$ |S| \ge |S+L| - |L| + 1 = |G| - |G_0| + 1 > \frac12\,|G|, $$
and it also follows that $G/L$ contains a subgroup isomorphic to
$(\Z/4\Z)\oplus(\Z/2\Z)$, which implies that $G$ itself contains such a
subgroup.
\end{proof}

Our next result shows that under the extra assumption $\kap(\CP(S))<|S|$,
the conclusion of Theorem \reft{main1} can be greatly simplified.
\begin{theorem}\label{t:simple}
Let $S$ be a proper subset of the finite abelian group $G$. If
$\kap(\CP(S))<|S|$, then
  $$ \kap(\CP(S)) = \min \{ |S+H|-|H|\colon H\le G,\, S+H\neq G \}. $$
\end{theorem}

Theorem \reft{simple} will be derived from Theorem \reft{main2} in Section
\refs{proofs}. Note that the assumption $\kap(\CP(S))<|S|$ of Theorem
\reft{simple} cannot be dropped: say, if $S$ is the non-zero coset of a
subgroup $H\le G$ of index $2$, then $\CP(S)$ is a complete bipartite graph
of connectivity $|G|/2$, while $|S+H|-|H|=0$ and $S+H\neq G$. We also notice
that, despite its simple and neat conclusion (and one which mirrors the
corresponding result for usual Cayley graphs), Theorem \reft{simple} gives no
way to determine whether $\kap(\CP(S))<|S|$ holds, and hence no way to find
the connectivity unless it is known to be smaller than $|S|$ a priori. Of
course, a necessary and sufficient condition for $\kap(\CP(S))<|S|$ to hold
follows readily from Theorem \reft{main2}.
\begin{corollary}\label{c:kappalarge}
If $S$ is a proper subset of the finite abelian group $G$, then in order for
$\kap(\CP(S))<|S|$ to hold it is necessary and sufficient that there is a
subgroup $K\in\cH_G(S)\cup\cL^\ast_G(S)$ with $|S+K|\le|S|+|K|-1$.
\end{corollary}

Observe that if $g$ is an element of $G$ with $2g\in S$, then $g$ is
a neighbor of itself in $\CP(S)$; consequently, the boundary of
$\{g\}$ contains $|S|-1$ elements so that $\kap(\CP(S))<|S|$. Hence
Theorem \reft{simple} implies the following corollary.

\begin{corollary}\label{c:2astS}
Let $S$ be a proper subset of the finite abelian group $G$. If
 $S\cap(2\ast G)\neq\est$, and in particular if $G$ has odd order and $S$ is
non-empty, then
  $$ \kap(\CP(S)) = \min \{ |S+H|-|H|\colon H\le G,\, S+H\neq G \}. $$
\end{corollary}

We conclude this section with two potentially useful lower-bound
estimates for $\kap(\CP(S))$.
\begin{corollary}\label{c:Sover2}
Let $S$ be a proper subset of the finite abelian group $G$. If $\CP(S)$
is connected, then in fact
  $$ \kap(\CP(S))\ge\frac12\,|S|. $$
\end{corollary}
Corollary \refc{Sover2} follows from Theorem \reft{simple} and the
observation that if $|S+H|-|H|=\kap(\CP(S))>0$ for a subgroup $H\le G$,
then $S$ intersects at least two cosets of $H$, so that $|S+H|\ge 2|H|$,
and therefore $|S+H|-|H|\ge\frac12\,|S+H|\ge\frac12\,|S|$.

\begin{corollary}\label{c:charG}
Let $S$ be a proper subset of the finite, non-trivial abelian group $G$, and
let $p$ denote the smallest order of a non-zero subgroup of $G$. If $\CP(S)$
is connected, then in fact
  $$ \kap(\CP(S)) \ge \min \{ |S|-1,\, p \}. $$
\end{corollary}
The proof is similar to that of the previous corollary: if
$\kap(\CP(S))<|S|-1$, then by Theorem \reft{simple} there exists a subgroup
$H\le G$ with $|S+H|-|H|=\kap(\CP(S))>0$; this subgroup is non-zero and hence
$|S+H|-|H|\ge|H|\ge p$.

\section{Auxiliary results}\label{s:aux}

In this section, we gather the tools needed for the proof of Theorems
\reft{main2} and \reft{simple}. This includes a simple consequence from
\refb{d1} or \refb{d2} (rephrased), a classical theorem of Kneser on
periodicity of sumsets, a result from \refb{l1}, which is a `dual' version of
a well-known structure theorem of Kemperman \refb{k}, and three original
lemmas.

Given a subgroup $H$ of the abelian group $G$, by $\pH$ we denote the
canonical homomorphism from $G$ onto $G/H$. Though the notation $\pH$ does
not specify the underlying group $G$, it is always implicit from the context
and no confusion will arise.

For a subset $S$ of the abelian group $G$, the (maximal) period of $S$ will
be denoted by $\pi(S)$; recall that this is the subgroup of $G$ defined by
  $$ \pi(S) := \{g\in G\colon S+g=S \}, $$
and that $S$ is called \emph{periodic} if $\pi(S)\neq\{0\}$ and
\emph{aperiodic} otherwise. Thus, $S$ is a union of $\pi(S)$-cosets, and
$\pi(S)$ lies above any subgroup $H\le G$ such that $S$ is a union of
$H$-cosets. Observe also that $\pi(S)=G$ if and only if either $S=\est$
or $S=G$, and that $\pH[\pi(G)](S)$ is an aperiodic subset of the group
$G/\pi(S)$.

\begin{theirproposition}[Grynkiewicz, \protect{\cite[(c.5)]{b:d1}; see also
                                \cite[Proposition 5.2]{b:d2}}]\label{l:david}
Let $A$ be a finite, non-empty subset of an abelian group. If
$|\pi(A\stm\{a\})|>2$ for some $a\in A$, then $|\pi(A\stm\{a'\})|=1$ for
every group element $a'\ne a$.
\end{theirproposition}

\begin{theirtheorem}[Kneser, \cite{b:kn1,b:kn2}; see also \refb{m}]
  \label{t:kneser}
Let $A$ and $B$ be finite, non-empty subsets of an abelian group $G$. If
  $$ |A+B| \le |A|+|B|-1, $$
then, letting $H:=\pi(A+B)$, we have
  $$ |A+B| = |A+H|+|B+H|-|H|. $$
\end{theirtheorem}

We now turn to the (somewhat involved) statement of \cite[Theorem 2]{b:l1};
the reader can consult the source for the explanations and comments.

By an arithmetic progression in the abelian group $G$ with difference
 $d\in G$, we mean a set of the form $\{g+d,g+2d\longc g+kd\}$, where $g$ is
an element of $G$ and $k$ is a positive integer. Thus, cosets of finite
cyclic subgroups (and in particular, singleton sets) are considered
arithmetic progressions, while the empty set is not. For finite subsets $A$
and $B$ of an abelian group and a group element $c$, we write
  $$ \nu_c(A,B) := |\{ (a,b)\in A\times B\colon c=a+b \}|; $$
that is, $\nu_c(A,B)$ is the number of representations of $c$ as a sum of an
element of $A$ and an element of $B$. Observe that $\nu_c(A,B)>0$ if and only
if $c\in A+B$. The smallest number of representations of an element of $A+B$
will be denoted by $\mu(A,B)$:
  $$ \mu(A,B) := \min \{ \nu_c(A,B) \colon c\in A+B \}. $$

Following Kemperman \refb{k}, we say that the pair $(A,B)$ of finite
subsets of the abelian group $G$ is \emph{elementary} if at least one of
the following conditions holds:
\begin{itemize}
\item[(I)]   $\min\{|A|,|B|\}=1$;
\item[(II)]  $A$ and $B$ are arithmetic progressions sharing a common
  difference, the order of which in $G$ is at least $|A|+|B|-1$;
\item[(III)] $A=g_1+(H_1\cup\{0\})$ and $B=g_2-(H_2\cup\{0\})$, where
  $g_1,g_2\in G$, and where $H_1$ and $H_2$ are non-empty subsets of a subgroup
  $H\le G$ such that $H=H_1\cup H_2\cup\{0\}$ is a partition of $H$;
  moreover, $c:=g_1+g_2$ is the only element of $A+B$ with $\nu_c(A,B)=1$;
\item[(IV)] $A=g_1+H_1$ and $B=g_2-H_2$, where $g_1,g_2\in G$,
  and where $H_1$ and $H_2$ are non-empty, aperiodic subsets of a subgroup $H\le G$
  such that $H=H_1\cup H_2$ is a partition of $H$; moreover, $\mu(A,B)\ge 2$.
\end{itemize}

We say that the pair $(A,B)$ of subsets of an abelian group satisfies
\emph{Kemperman's condition} if either $A+B$ is aperiodic or $\mu(A,B)=1$
holds.

\begin{theirtheorem}[\protect{Lev, \cite[Theorem 2]{b:l1}}]\label{t:dual}
Let $A$ and $B$ be finite, non-empty subsets of the abelian group $G$. A
necessary and sufficient condition for $(A,B)$ to satisfy both
  $$ |A+B| \le |A| + |B| - 1 $$
and Kemperman's condition is that there exist non-empty subsets
 $A_0\seq A$ and $B_0\seq B$ and a finite, proper subgroup $F<G$ such that
\begin{itemize}
\item[(i)]   each of $A_0$ and $B_0$ is contained in an $F$-coset,
  $|A_0+B_0|=|A_0|+|B_0|-1$, and the pair $(A_0,B_0)$ satisfies Kemperman's
  condition;
\item[(ii)]  each of $A\stm A_0$ and $B\stm B_0$ is a (possibly empty) union
  of $F$-cosets;
\item[(iii)] the pair $(\pH[F](A),\pH[F](B))$ is elementary; moreover, either
  $F$ is trivial, or $\pH[F](A_0)+\pH[F](B_0)$ has a unique representation as
  a sum of an element of $\pH[F](A)$ and an element of $\pH[F](B)$.
\end{itemize}
\end{theirtheorem}

\begin{lemma}\label{l:directsum}
Let $L\le G_0\le G$ be finite abelian groups. If $G_0/L$ is a cyclic
$2$-group and $2\ast(G/L)$ is a proper subgroup of $G_0/L$, then
$\exp(G_0/L)=\exp(G/L)$.
\end{lemma}

\begin{proof}
Write $|G_0/L|=2^k$ so that $k$ is a positive integer. Since $|2\ast(G/L)|$
is a proper divisor of $2^k$, we have $2^{k-1}g=0$ for every
 $g\in 2\ast(G/L)$. Equivalently, $2^kg\in L$ for every $g\in G$, whence
$\exp(G/L)\le 2^k=\exp(G_0/L)$. The inverse estimate
$\exp(G_0/L)\le\exp(G/L)$ is trivial.
\end{proof}

The following lemma is similar in flavor to a lemma used by Kneser to prove
Theorem \reft{kneser}; cf. \cite{b:kn2,b:k}.

\begin{lemma}\label{l:geom}
Suppose that $S$ is a finite subset, and that $H$ and $L$ are finite
subgroups of the abelian group $G$ satisfying $|L|\le|H|$ and $S+H\ne S+H+L$.
Let $I:=H\cap L$. If
  $$ \max \{ |S+H| - |H|, |S+L| - |L| \} \le |S+I| - |I|, $$
then in fact
  $$ |S+H| - |H| = |S+L| - |L| = |S+I| - |I|; $$
moreover, there exists $g\in G$ such that $(S+I)\stm(g+H+L)$ is a (possibly
empty) union of $(H+L)$-cosets, and one of the following holds:
\begin{itemize}
\item[(i)]  $(S+I)\cap(g+H+L)=g+I$;
\item[(ii)] $(S+I)\cap(g+H+L)=(g+H+L)\stm(g+(H\cup L))$ and $|H|=|L|$.
\end{itemize}
\end{lemma}

\begin{proof}
Factoring by $I$, we assume without loss of generality that $I=\{0\}$. Since
$S+H\neq S+H+L$, there exists $s_0\in S$ with $s_0+L\nsubseteq S+H$, and we
let $S_0:=S\cap(s_0+H+L)$. It is instructive to visualize the coset $s_0+H+L$
as the grid formed by $|L|$ horizontal lines (corresponding to the $H$-cosets
contained in $s_0+H+L$) and $|H|$ vertical lines (corresponding to the
$L$-cosets contained in $s_0+H+L$). The intersection points of these two
families of lines correspond to the elements of $s_0+H+L$, and the condition
$s_0+L\nsubseteq S+H$ implies that there is a horizontal line free of
elements of $S$.

Let $h:=\pH[L](S_0)$ (the number of vertical lines that intersect $S_0$) and
$l:=\pH(S_0)$ (the number of horizontal lines that intersect $S_0$); thus,
$1\le h\le |H|$ and $1\le l<|L|$. We also have, in view of the hypotheses,
\begin{equation}\label{e:geom0}
  (|H|-h)l \le |(S_0+H)\stm S_0| \le |(S+H)\stm S| \le |H| - 1,
\end{equation}
whence
\begin{equation}\label{e:geom1}
  (|H|-h)(l-1) \le h-1,
\end{equation}
and similarly,
\begin{equation}\label{e:geom2}
  (|L|-l)(h-1) \le l-1.
\end{equation}

To begin with, suppose that $l=1$, and hence $h=1$ by \refe{geom2}. In this
case, $|S_0|=1$, whence $S\cap(s_0+H+L)=\{s_0\}$. Furthermore, \refe{geom0}
yields $(S_0+H)\stm S_0=(S+H)\stm S$, and likewise we have
 $(S_0+L)\stm S_0=(S+L)\stm S$. This shows that
\begin{equation}\label{e:geom2a}
  |S+H|-|S|=|H|-1, \quad |S+L|-|S|=|L|-1,
\end{equation}
and $S\stm S_0$ is a union of $(H+L)$-cosets, thus establishing the assertion
(with $g=s_0$) in the case $l=1$. So we assume $l>1$ below.

Observe that \refe{geom1} and \refe{geom2} imply
  $$ l-1 \ge (|L|-l)(h-1) \ge (|L|-l)(|H|-h)(l-1), $$
whence it follows from $l>1$ that
\begin{equation}\label{e:geom3}
  (|L|-l)(|H|-h) \le 1.
\end{equation}
If $|H|=h$, then \refe{geom2} gives
  $$ l-1 \ge (|H|-1)(|L|-l) \ge (|L|-1)(|L|-l) \ge 2|L|-l-2 \ge l, $$
which is wrong. Therefore $|H|>h$. Thus we deduce from \refe{geom3}
and $l<|L|$ that $h=|H|-1$ and $l=|L|-1$, whence \refe{geom2} gives
$|H|=|L|$. Consequently, \refe{geom0} yields $(S_0+H)\stm
S_0=(S+H)\stm S$, and similarly $(S_0+L)\stm S_0=(S+L)\stm S$, which
(as above) proves \refe{geom2a} and shows that $S\stm S_0$ is a
union of $(H+L)$-cosets. Furthermore, $S+H$ misses exactly one
$H$-coset in $s_0+H+L$, and $S+L$ misses exactly one $L$-coset in
$s_0+H+L$. Let $g\in s_0+H+L$ be the common element of these two
cosets, so that $S_0+H=(s_0+H+L)\stm(g+H)$ and
$S_0+L=(s_0+H+L)\stm(g+L)$. Then
  $$ S_0 \seq (s_0+H+L)\stm(g+(H\cup L)) = (g+H+L)\stm(g+(H\cup L)), $$
and thus
  $$ |L|-1 = |H|-1 \ge |(S+H)\stm S| = |(S_0+H)\stm S_0|
                                                 = (|L|-1)|H| - |S_0|, $$
so that
  $$ |S_0| \ge (|H|-1)(|L|-1) = |(g+H+L)\stm(g+(H\cup L))|.$$
Hence, in fact $S_0=(g+H+L)\stm(g+(H\cup L))$, completing the proof.
\end{proof}

\begin{lemma}\label{l:LcapH}
Let $G$ be a finite abelian group, and suppose that the proper
subset $S\subset G$, the subgroups $L\le G_0\le G$, and the element
$g_0\in G_0$ satisfy conditions (L1)--(L4) in the definition of
$\cL^\ast_G(S)$. Suppose, moreover, that $|S+L|-|S|\le|L|-1$. If $H$
is a subgroup of $G$ with $|S+H|-|S|\le|H|-1$ and $S+H\ne G$, then
$H$ is actually a subgroup of $G_0$.
\end{lemma}

\begin{proof}
Suppose for a contradiction that $H\nleq G_0$ and fix $h\in H\stm G_0$. For
each $g\in G_0$, we have $g+h\in G\stm G_0\seq S+L$, whence $g\in S+H+L$.
Hence $G_0\seq S+H+L$, and since, on the other hand, we have
 $G\stm G_0\seq S+L\seq S+H+L$, we conclude that
\begin{equation}\label{e:LcapH1}
  S+H+L = G.
\end{equation}
In view of $S+H\ne G$, this leads to $L\nleq H$, and we let $I:=H\cap L$.
Thus $I$ is a proper subgroup of $L$.

Write $n:=|G_0/L|$ so that $G_0$ consists of $n\ge 4$ cosets of $L$, of which
$n-1$ are free of elements of $S$. Let $\{g_i\colon 0\le i\le n-1\}$ be a
system of representatives of these $n$ cosets.

Fix $i\in[1,n-1]$. Since $H\nleq G_0$ and $g_i\in G_0$, we have
$g_i+H\nsubseteq G_0$, whence $(G\stm G_0)\cap(g_i+H)\neq\est$; as
$G\stm G_0\seq S+L$, this yields $S\cap(g_i+H+L)\neq\est$. On the
other hand, from $g_i+L\seq G_0\stm(g_0+L)$ it follows that
$(S+L)\cap (g_i+L)=\est$. Therefore,
\begin{equation}\label{e:LcapH2}
  0 < |(S+I)\cap(g_i+H+L)| < |H+L|; \quad i\in[1,n-1].
\end{equation}

In view of \refe{LcapH1} and the hypotheses $S+H\neq G$, we have
 $S+H\ne S+H+L$ and $S+L\ne S+H+L$. Also, our assumptions imply
  $$ \max \{ |S+H|-|H|, |S+L|-|L| \} < |S| \le |S+I|, $$
and since both the left and right hand side are divisible by $|I|$, we
actually have
  $$ \max \{ |S+H|-|H|, |S+L|-|L| \} \le |S+I|-|I|. $$
Thus we can apply Lemma \refl{geom}. Choose $g\in G$ such that
$(S+I)\stm(g+H+L)$ is a union of $(H+L)$-cosets. Then it follows from
\refe{LcapH2} that
\begin{equation}\label{e:LcapH3}
  g_i+H+L = g+H+L; \quad i\in[1,n-1],
\end{equation}
and consequently $G_0\stm(g_0+L)\seq g+H+L$. Hence $n\ge 4$ implies
 $G_0\le H+L$ and $g\in H+L$. Thus, since $S\cap(g_0+L)$ is not contained in a
coset of a proper subgroup of $L$, and in particular in a coset of $I$, we
conclude that
  $$ |(S+I)\cap(g+H+L)| = |(S+I)\cap(g_0+L)| \ge 2|I|. $$
This shows that Lemma \refl{geom}~(i) fails. On the other hand, \refe{LcapH3}
gives $g_i+L\seq g+H+L$, and hence $g+H+L$ contains at least $n-1\ge 3$
cosets of $L$, all free of elements of $S+I$. Thus Lemma \refl{geom}~(ii)
fails too, a contradiction.
\end{proof}

\section{Proofs of Theorems \reft{main2} and \reft{simple}}\label{s:proofs}

Our starting point is the observation that if $S$ is a subset of the finite
abelian group $G$ such that $\CP(S)$ is not complete, then
  $$ \kap(\CP(S)) = \min \{ |(S-A)\stm A|
                         \colon \est\ne A\seq G,\, (S-A)\cup A\ne G \}. $$

For the following proposition, the reader may need to recall the notion of a
fragment, introduced in Section \refs{intro} after the statement of Theorem
\reft{main1}.
\begin{proposition}\label{p:lessS}
Let $S$ be a subset of the finite abelian group $G$, and suppose that
$\kap(\CP(S))<|S|$. If $A$ is a fragment of $\CP(S)$, then, writing
$H:=\pi(S-A)$, we have
\begin{align}
  A                      &\seq S-A,        \label{e:-t}\\
  A+H                    &=    A,          \label{e:t+h}\\
  \kap(\CP(S))           &=    |S+H|-|H|,  \label{e:per}\\
\intertext{and}
  \kap(\CP[G/H]{\pH(S)}) &=    |\pH(S)|-1. \label{e:g/h}
\end{align}
\end{proposition}

\begin{proof}
Fix $a\in A$. Since $a$ has $|S|$ neighbors, all lying in $S-A$, and since
$|(S-A)\stm A|=\kap(\CP(S))<|S|$ by the assumptions, it follows that $a$ has
a neighbor in $A$; in other words, there is $a'\in A$ with $a+a'\in S$.
Consequently, $a\in S-A$, and \refe{-t} follows.

By \refe{-t} we have
  $$ (S-(A+H))\cup(A+H) = S-A+H = S-A \neq G, $$
and obviously,
  $$ |(S-(A+H)) \stm (A+H)| = |(S-A) \stm (A+H)| \le |(S-A) \stm A|. $$
Since $A$ is a fragment, we conclude that in fact
 $|(S-A)\stm(A+H)|=|(S-A)\stm A|$ holds, which gives \refe{t+h}.

By \refe{-t} and the assumptions, we have
  $$ |S-A| = |(S-A)\stm A| + |A| = \kap(\CP(S)) + |A| \le |S|+|A|-1. $$
Hence it follows from Theorem \reft{kneser} and \refe{t+h} that
\begin{equation}\label{e:loc1}
  |S-A| = |S+H| + |A+H| - |H| = |S+H| + |A| - |H|.
\end{equation}
Thus
  $$ \kap(\CP(S)) = |(S-A)\stm A| = |S-A|-|A| = |S+H|-|H|, $$
yielding (\ref{e:per}).

Finally, we establish \refe{g/h}. The neighborhood of $\pH(A)$ in the
graph $\CP[G/H](\pH(S))$ is $\pH(S)-\pH(A)=\pH(S-A)$, and  it follows in
view of \refe{-t} that
  $$ \pH(S-A) \cup \pH(A) = \pH(S-A) \neq G/H.$$
 Consequently, the set $\pH(S-A)\stm\pH(A)$ is a vertex cut in
$\CP[G/H](\pH(S))$, whence using \refe{-t}, \refe{t+h}, and \refe{loc1}
we obtain
\begin{multline*}
  \kap(\CP[G/H](\pH(S))) \le |\pH(S-A)\stm\pH(A)|=|\phi_H(S-A)|-|\phi_H(A)| \\
                        = (|S-A|-|A|)/|H| = |S+H|/|H| - 1 = |\pH(S)| - 1.
\end{multline*}
To prove the inverse estimate, notice that the graph $\CP[G/H](\pH(S))$ is
not complete (we saw above that it has vertex cuts) and choose $A'\seq G$
such that $\pH(A')$ is a fragment of this graph. Replacing $A'$ with $A'+H$,
we can assume without loss of generality that $A'+H=A'$. Since
  $$ \pH((S-A')\cup A') = (\pH(S) - \pH(A')) \cup \pH(A') \ne G/H, $$
we have $(S-A')\cup A'\neq G$. Hence in view of \refe{per} it follows that
\begin{align*}
  \kap(\CP[G/H](\pH(S)))
    &=   |(\pH(S)-\pH(A'))\stm \pH(A')| \\
    &=   |\pH(S-A')\stm\pH(A')| \\
    &=   |(S-A')\stm A'|/|H| \\
    &\ge |\kap(\CP(S))|/|H| \\
    &=   |\pH(S)|-1,
\end{align*}
as desired.
\end{proof}

For a subset $S$ of a finite abelian group $G$, write
  $$ \lam^\ast_G(S) := \min \{ |S+L|-|L| \colon L\in\cL^\ast_G(S) \}. $$
Clearly, we have $\lam^\ast_G(S)\ge\lam_G(S)$.

\begin{lemma}\label{l:Shift}
Let $S$ be a proper subset of the finite abelian group $G$. If $g\in G$,
then $\cH_G(S-2g)=\cH_G(S),\ \cL^\ast_G(S-2g)=\cL^\ast_G(S)$, and
$\CP(S-2g)$ is isomorphic to $\CP(S)$; consequently,
\begin{gather*}
  \eta_G(S-2g)=\eta_G(S),\ \lam^\ast_G(S-2g)=\lam^\ast_G(S), \\
\intertext{and}
  \kap(\CP(S-2g))=\kap(\CP(S)).
\end{gather*}
\end{lemma}

\begin{proof}
The isomorphism between $\CP(S-2g)$ and $\CP(S)$ is established by mapping
every group element $x$ to $x-g$, and the equality $\cH_G(S-2g)=\cH_G(S)$ is
immediate from the observation that $S+2\ast G-2g=S+2\ast G$. To show that
$\cL^\ast_G(S-2g)=\cL^\ast_G(S)$, suppose that $L\in\cL^\ast_G(S)$ and let
$G_0\le G$ (lying above $L$) and $g_0\in G_0$ be as in (L1)--(L4). By (L2) we
have $2g\in G_0$. Consequently, $(G\stm G_0)-2g=G\stm G_0$, and hence it
follows from (L4) that
  $$ S-2g+L = (G\stm G_0) \cup (g_0-2g+L). $$
Furthermore, since $\pH[L](g_0)$ is a generator of the cyclic $2$-group
$G_0/L$, so is $\pH[L](g_0-2g)$; that is, $\<g_0-2g\>+L=G_0$. This shows
that $L\in\cL^\ast_G(S-2g)$, and hence
$\cL^\ast_G(S)\seq\cL^\ast_G(S-2g)$. By symmetry, we also have
$\cL^\ast_G(S-2g)\seq\cL^\ast_G(S)$, implying the assertion.
\end{proof}

We now pass to our last lemma, which will take us most of the way towards
the proof of Theorem \reft{main2}; the reader may compare the statement
of this lemma with that of Theorem \reft{main1}.
\begin{lemma}\label{l:main}
If $S$ is a proper subset of the finite abelian group $G$, then
\begin{equation}\label{e:THEresult}
  \kap(\CP(S)) = \min \{ \eta_G(S), \lam^\ast_G(S), |S| \}.
\end{equation}
\end{lemma}

\begin{proof}
Since each of $\eta_G(S),\,\lam^\ast_G(S)$, and $|S|$ is an upper bound for
$\kap(\CP(S))$, it suffices to show that $\kap(\CP(S))$ is greater than or
equal to one of these quantities. Thus we can assume that
 $\kap(\CP(S))\le |S|-1\le |G|-2$. Hence $S\ne\est$ and $\CP(S)$ is not
complete.

It is not difficult to see that the assertion holds true if $|G|\le 2$; we
leave verification to the reader. The case $|S|=1$ is also easy to establish
as follows. Suppose that $|G|>2$ and $S=\{s\}$, where $s$ is an element of
$G$. If $\<s\>\ne G$, then $\<s\>\in\cH_G(S)$ and $|S+\<s\>|-|\<s\>|=0$,
implying $\kap(\CP(S))=\eta_G(S)=0$. Next, if $G$ is not a $2$-group, then
there exists an element $g\in G$ which is an odd multiple of $s$ and such
that the subgroup $\<g\>$ is proper; in this case $g\in(S+2\ast G)\cap\<g\>$
showing that $\<g\>\in\cH_G(S)$ and leading to $\kap(\CP(S))=\eta_G(S)=0$, as
above. In both cases the proof is complete, so we assume that $\<s\>=G$ is a
$2$-group. Since $|G|>2$, in this case we have $\{0\}\in\cL^\ast_G(S)$ (take
$G_0=G$ and $g_0=s$ in (L1)--(L4)) and $|S+\{0\}|-|\{0\}|=0$, whence
$\kap(\CP(S))=\lam^\ast_G(S)=0$.

Having finished with the cases $|S|=1$ and $|G|\le 2$, we proceed by
induction on $|G|$, assuming that $\kap(\CP(S))\le|S|-1$. Choose $A\seq G$
such that $A$ is a fragment of $\CP(S)$ and fix arbitrarily $a\in A$. In view
of Lemma \refl{Shift}, and since the set $A-a$ is a fragment of the graph
$\CP(S-2a)$, by passing from $S$ to $S-2a$, and from $A$ to $A-a$, we ensure
that
\begin{equation}\label{e:0inT}
  0 \in A.
\end{equation}
Also, by Proposition \refp{lessS} we have $A\seq S-A\ne G$.

If each of $S$ and $A$ is contained in a coset of a proper subgroup $K<G$,
then from $A\seq S-A$ and \refe{0inT} it follows that in fact $S$ and $A$ are
contained in $K$, whence $K\in\cH_G(S)$; furthermore, $|S+K|-|K|=0$, showing
that $\kap(\CP(S))=\eta_G(S)=0$. Accordingly, we assume for the rest of the
proof that for any proper subgroup of $G$, at least one of the sets $S$ and
$A$ is not contained in a coset of this subgroup.

Let $H:=\pi(S-A)$. We distinguish two major cases according to whether or not
$H$ is trivial.

\subsection*{Case 1:} $H$ is non-trivial. Applying the induction
hypothesis to $\CP[G/H](\pH(S))$ and using \refe{g/h}, we conclude
that either $\eta_{G/H}(\pH(S))=|\pH(S)|-1$ or
$\lam^\ast_{G/H}(\pH(S))=|\pH(S)|-1$, giving two subcases.

\subsubsection*{Subcase 1.1.}
Assume first that $\eta_{G/H}(\pH(S))=|\pH(S)|-1$, and hence that there
exists a subgroup $H'\le G$, lying above $H$, such that
$H'/H\in\cH_{G/H}(\pH(S))$ and
  $$ |\pH(S)+H'/H| - |H'/H| = \eta_{G/H}(\pH(S)) = |\pH(S)|-1. $$
The former easily implies that $H'\in\cH_G(S)$, while the latter, in
conjunction with \refe{per}, implies that
  $$ |S+H'| - |H'| = |S+H| - |H| = \kap(\CP(S)). $$
This shows that $\kap(\CP(S))\ge\eta_G(S)$, whence in fact
$\kap(\CP(S))=\eta_G(S)$.

\subsubsection*{Subcase 1.2.}
Assume now that $\lam^\ast_{G/H}(\pH(S))=|\pH(S)|-1$, and let $L\le G$ be
a subgroup, lying above $H$, such that $L/H\in\cL^\ast_{G/H}(\pH(S))$ and
  $$ |\pH(S)+L/H| - |L/H| = \lam^\ast_{G/H}(\pH(S)) = |\pH(S)|-1. $$
In view of \refe{per} and the assumptions, the last equality yields
\begin{equation}\label{e:yikess}
  |S+L|-|L| = |S+H|-|H| = \kap(\CP(S)) \le |S|-1.
\end{equation}

Since $L/H\in\cL^\ast_{G/H}(\pH(S))$, we can find a subgroup $G_0\le G$,
lying above $L$, and an element $g_0\in G_0\stm L$, so that $G/G_0$ is an
elementary abelian $2$-group, $G_0/L$ is a cyclic $2$-group of order at least
$4$ generated by $\pH[L](g_0)$, and $S+L=(G\stm G_0)\cup(g_0+L)$. Without
loss of generality, we can assume that $g_0\in S$.

If $S_0:=S\cap(g_0+L)$ is not contained in a coset of a proper subgroup of
$L$, then $L\in\cL^\ast_G(S)$, and hence it follows in view of \refe{yikess}
that $\kap(\CP(S))=\lam^\ast_G(S)$. Therefore we assume that there exists a
proper subgroup $R<L$ such that $S_0$ is contained in an $R$-coset, and we
choose $R$ to be minimal subject to this property; thus, $S_0=S\cap(g_0+R)$
and $\<(S-g_0)\cap L\>=R$.

Since $S_0$ is contained in an $R$-coset, from \refe{yikess} we obtain
  $$ |(S\stm S_0)+L| - |S\stm S_0|
                               = |S+L|-|L| - |S| + |S_0| < |S_0| \le |R|. $$
Hence every $R$-coset in
 $G\stm G_0=(S\stm S_0)+L$ contains at least one element of $S$; that is,
\begin{equation}\label{e:S+R=}
   S+R = (G\stm G_0)\cup(g_0+R).
\end{equation}
Consequently, using \refe{yikess} once again, we obtain
\begin{equation}\label{e:S+R-R}
  |S+R|-|R| = |G\stm G_0| = |S+L|-|L| = \kap(\CP(S)).
\end{equation}

Applying the previously completed singleton case to the set
 $\pH[R](S_0)\seq G_0/R$, we get two further subcases.

\subsubsection*{Subcase 1.2.1.}

Suppose that $\kap(\CP[G_0/R](\pH[R](S_0)))=\eta_{G_0/R}(\pH[R](S_0))$.
Choose a subgroup $R'\le G_0$, lying above $R$, such that
$R'/R\in\cH_{G_0/R}(\pH[R](S_0))$. Since $R\le R'\le G_0$, it follows in view
of \refe{S+R=} and \refe{S+R-R} that
  $$ |S+R'|-|R'| = |S+R|-|R| = \kap(\CP(S)). $$
Thus, since $R'\in\cH_{G_0}(S_0)\seq\cH_G(S)$, we conclude that
$\kap(\CP(S))=\eta_G(S)$.

\subsubsection*{Subcase 1.2.2.}
Assume now that $\kap(\CP[G_0/R](\pH[R](S_0)))\ne\eta_{G_0/R}(\pH[R](S_0))$.
As $|G_0/R|\ge|G_0/L|\ge 4$, from the singleton case analysis at the
beginning of the proof it follows that $G_0/R$ is a cyclic $2$-group
generated by $\pH[R](S_0)=\{\pH[R](g_0)\}$.

If $R\in\cH_G(S)$, then it follows in view of \refe{S+R-R} that
$\kap(\CP(S))=\eta_G(S)$; therefore, we assume that $R\notin\cH_G(S)$. Hence
in view of $S+R\seq S+L\ne G$ we infer that $2\ast(G/R)\cap\pH[R](S)=\est$.
Consequently, since \refe{S+R=} implies that $\pH[R](S)$ contains
$(G/R)\stm(G_0/R)$ as a proper subset, we have $2\ast(G/R)\lneqq G_0/R$.

Applying Lemma \refl{directsum}, we conclude that $\exp(G_0/R)=\exp(G/R)$.
Thus \refe{S+R=}, the remark at the beginning of the present subcase, and the
above-made observation that $G/G_0$ is an elementary $2$-group show that
$R\in\cL^\ast_G(S)$, whence \refe{S+R-R} yields
$\kap(\CP(S))=\lam^\ast_G(S)$.

\subsection*{Case 2:} $H$ is trivial. Thus by \refe{per} we have
$\kap(\CP(S))=|S|-1$, and therefore \refe{-t} gives
  $$ |S-A|-|A| = |(S-A)\stm A| = \kap(\CP(S)) = |S|-1. $$
Applying Theorem \reft{dual} to the pair $(S,-A)$, we find a subgroup $F<G$
such that conclusions (i)--(iii) of Theorem \ref{t:dual} hold true; in
particular, $(\pH[F](S),-\pH[F](A))$ is an elementary pair in $G/F$ of one of
the types (I)--(IV), and $|S+F|\le |S|+|F|-1$. By the last inequality, we
have
  $$ |S+F|-|F| \le |S|-1 = \kap(\CP(S)). $$
Hence, if $F\in\cH_G(S)$, then $\kap(\CP(S))=\eta_G(S)$; consequently, we
assume that
\begin{equation}\label{e:FnotincH}
  F\notin\cH_G(S).
\end{equation}

Observe that if $\pH[F](S)=G/F$, then $F$ is non-zero, whence by Theorem
\reft{dual} (iii) we have $|\pH[F](A)|=1$. Thus, if $(\pH[F](S),-\pH[F](A))$
is not of type (I), then
\begin{equation}\label{e:technicality-condition}
  S+F \ne G.
\end{equation}
We proceed by cases corresponding to the type of the pair
$(\pH[F](S),-\pH[F](A))$.

\subsubsection*{Subcase 2.1.}
Suppose that $(\pH[F](S),-\pH[F](A))$ is of type (IV). In this case, we have
$\mu(\pH[F](S),-\pH[F](A))\ge 2$, whence it follows by
Theorem~\reft{dual}~(iii) that $F$ is trivial. Hence $(S,-A)$ is an
elementary pair of type (IV). Thus, since $S$ and $A$ are not both contained
in a coset of the same proper subgroup, it follows that $A=g+(G\stm S)$ for
some $g\in G$, implying $-g\notin S-A$. Therefore \refe{-t} yields
 $-g\notin g+(G\stm S)$ and thus $-2g\in S$; consequently,
$\{0\}=F\in\cH_G(S)$, contradicting \refe{FnotincH}.

\subsubsection*{Subcase 2.2.}
Suppose that $(\pH[F](S),-\pH[F](A))$ is of type (III), but not of type (I).
Then, since $S$ and $A$ are not both contained in a coset of the same proper
subgroup and since $S-A\neq G$, it follows that $F$ is non-zero, that
  $$ \pH[F](S) = \pH[F](g_1) + (H_1\cup\{0\})
                            ,\ -\pH[F](A) = \pH[F](g_2) - (H_2\cup\{0\}) $$
for some $g_1,g_2\in G$, where $H_1\cup H_2\cup\{0\}$ is a partition
of $G/F$, and that $g_1+g_2+F$ has a non-empty intersection with
$S-A$, while every $F$-coset, other than $g_1+g_2+F$, is contained
in $S-A$; moreover, from $\pi(S-A)=\{0\}$ we derive that
\begin{equation}\label{e:r_0-not-a-subset}
  g_1+g_2+F\nsubseteq S-A.
\end{equation}
By Theorem \reft{dual}, all $F$-cosets corresponding to
  $$ (-\pH[F](A)) \stm \{\pH[F](g_2)\} = \pH[F](g_2)-H_2,$$
are contained in $-A$. Hence, if
  $$ -\pH[F](g_1+g_2) \in \pH[F](g_2)-H_2, $$
then $-g_1-g_2+F\seq -A$, and it follows in view of \refe{-t} that
$g_1+g_2+F\seq A\seq S-A$, contradicting \refe{r_0-not-a-subset}. Therefore,
assume instead that $-\pH[F](g_1+g_2)\notin\pH[F](g_2)-H_2$, so that
  $\pH[F](g_1+2g_2)\in H_1\cup\{0\}$.
Then $2\pH[F](g_1+g_2)\in\pH[F](g_1)+(H_1\cup\{0\})=\pH[F](S)$,
whence by \refe{technicality-condition} we have $F\in\cH_G(S)$,
contradicting \refe{FnotincH}.

\subsubsection*{Subcase 2.3.}
Suppose that $(\pH[F](S),-\pH[F](A))$ is of type (II), but not of type (I).
Letting $u:=|\pH[F](S)|$ and $v:=|\pH[F](A)|$, and choosing
 $s_0\in S,\,a_0\in A$, and
$d\in G\stm\{0\}$ appropriately, we write
\begin{align*}
   \pH[F](S) &= \{ \pH[F](s_0), \pH[F](s_0)+\pH[F](d),
                                  \ldots, \pH[F](s_0)+(u-1)\pH[F](d) \} \\
\intertext{and}
  -\pH[F](A) &= \{ \pH[F](a_0), \pH[F](a_0)+\pH[F](d),
                                  \ldots, \pH[F](a_0)+(v-1)\pH[F](d) \}.
\end{align*}
Since $(\pH[F](S),-\pH[F](A))$ is not of type (I), we have $u,\,v\ge 2$.
Next, it follows from \refe{-t} that
  $$ -\pH[F](a_0) = \pH[F](s_0)+\pH[F](a_0) + r\pH[F](d), $$
and therefore $\pH[F](s_0)=-2\pH[F](a_0)-r\pH[F](d)$, for some integer $r$.
Thus either $\pH[F](s_0)$ (if $r$ is even) or $\pH[F](s_0)+\pH[F](d)$ (if $r$
is odd) belongs to $2\ast(G/F)$. In either case, in view of $u\ge 2$ we have
$\pH[F](S)\cap(2\ast(G/F))\neq\est$, which by \refe{technicality-condition}
leads to $F\in\cH_G(S)$, contradicting \refe{FnotincH}.

\subsubsection*{Subcase 2.4.}
Finally, suppose that $(\pH[F](S),-\pH[F](A))$ is of type (I); that is,
either $|\pH[F](S)|=1$ or $|\pH[F](A)|=1$ holds.

Suppose first that $|\pH[F](S)|=1$. In this case, $F$ is non-zero (as
$|S|>1$) and $S+F\neq G$ (as $F$ is a \emph{proper} subgroup); moreover, from
\refe{-t} we obtain
\begin{equation}\label{e:cackle-cackle}
  \pH[F](S)-\pH[F](A) = \pH[F](A).
\end{equation}
By Theorem \reft{dual}, we can write $A=A_1\cup A_0$, where $A_1$ is a union
of $F$-cosets and $A_0$ is a non-empty subset of an $F$-coset disjoint from
$A_1$. If $\pH[F](S)-\pH[F](A_0)\seq\pH[F](A_1)$, then
 $S-A_0+F\seq A_1+F=A_1\seq S-A$, whence $S-A=(S-A_1)\cup (S-A_0)$
is a union of $F$-cosets, contradicting the assumption that $S-A$ is
aperiodic. Therefore \refe{cackle-cackle} gives
$\pH[F](S)-\pH[F](A_0)=\pH[F](A_0)$, which together with $S+F\ne G$ implies
$F\in\cH_G(S)$, contradicting \refe{FnotincH}. So we assume for the remainder
of the proof that $|\pH[F](S)|>|\pH[F](A)|=1$, and consequently in view of
\refe{0inT} that $A\seq F$.

Thus from \refe{-t} we derive that $0\in\pH[F](S)$, and it follows in view of
\refe{FnotincH} that $S+F=G$. Hence $F$ is nontrivial, and Theorem
\reft{dual} shows that there exists $s_0\in S$ such that $S=(G\stm
(s_0+F))\cup S_0$, where $S_0\subset s_0+F$.

If there exists $g\in G$ with $\pH[F](g)\ne -\pH[F](g)+\pH[F](s_0)$, then it
follows in view of $\pH[F](S)=G/F$ that
 $\pH[F](g)\in -\pH[F](g)+\pH[F](S\stm S_0)$, whence
  $$ g\in -g+(S\stm S_0)+F \seq -g+S; $$
consequently, $\{0\}\in\cH_G(S)$ and $\kap(\CP(S))=\eta_G(S)$. Therefore we
assume that $\pH[F](g)=-\pH[F](g)+\pH[F](s_0)$ for all $g\in G$. Hence
$2\ast(G/F)=\{\pH[F](s_0)\}$, which implies that $G/F$ is an elementary
$2$-group and that $\pH[F](s_0)=0$; consequently, $S_0=S\cap F$.

From $A\seq F$ and \refe{-t}, it follows that $A\seq(S-A)\cap F=S_0-A$, and
since $S-A\ne G$ and $S+F=G$ we have $S_0-A\ne F$. Consequently, Theorem
\reft{dual} (i) yields
\begin{equation}\label{e:statue}
  \kap(\CP[F](S_0)) \le |(S_0-A)\stm A| = |S_0-A|-|A| \le |S_0|-1.
\end{equation}
Since $S_0$ is a proper subset of $F$, it follows in view of \refe{statue}
that $\kap(\CP[F](S_0))\le |F|-2$, whence $\CP[F](S_0)$ is not complete. Let
$A'\seq F$ be a fragment of $\CP[F](S_0)$. By \refe{-t} and \ref{e:statue},
we have $A'\seq S_0-A'\ne F$, and consequently $A'\seq S-A'\ne G$. Hence from
\refe{statue} and $S\stm S_0=G\stm F$ we obtain
\begin{multline*}
  |S|-1 = \kap(\CP(S)) \le |(S-A')\stm A'|
                                          \le |G\stm F|+|(S_0-A')\stm A'| \\
        = |S\stm S_0|+\kap(\CP[F](S_0)) \le |S|-1,
\end{multline*}
implying $\kap(\CP[F](S_0))=|S_0|-1$ and
  $$ \kap(\CP(S)) = |S\stm S_0|+\kap(\CP[F](S_0)). $$
Consequently, if $F'\leq F$ has the property that
$\kap(\CP[F](S_0))=|S_0+F'|-|F'|$, then
\begin{equation}\label{e:statue2}
  \kap(\CP(S))=|S+F'|-|F'|.
\end{equation}

With \refe{statue} in mind, we apply the induction hypothesis to the
graph $\CP[F](S_0)$. If $\kap(\CP[F](S_0))=\eta_F(S_0)$, then by
\refe{statue2} any subgroup $F'\in\cH_F(S_0)\seq\cH_G(S)$ with
$\kap(\CP[F](S_0))=|S_0+F'|-|F'|$ satisfies
$\kap(\CP(S))=|S+F'|-|F'|$, whence $\kap(\CP(S))=\eta_G(S)$.
Therefore we assume instead that
$\kap(\CP[F](S_0))=\lam^\ast_F(S_0)$.

Choose $L\in\cL^\ast_F(S_0)$ with $\lam^\ast_F(S_0)=|S_0+L|-|L|$, and let
$G_0$ and $g_0\in G_0$ be as in (L1)--(L4), with $F$ playing the role of $G$.
Then it follows in view of \refe{statue2} that
\begin{equation}\label{e:k}
  \kap(\CP(S)) = |S+L|-|L|.
\end{equation}
If $\pH[L](S)\cap2\ast(G/L)\ne\est$, then $L\in\cH_G(S)$, whence \refe{k}
yields $\kap(\CP(S))=\eta_G(S)$. Therefore we assume that
\begin{equation}\label{e:hren1}
  \pH[L](S) \cap 2\ast (G/L)=\est
\end{equation}
and we proceed to show that $L\in\cL^\ast_G(S)$; in view of \refe{k}, this
will complete the proof.

Since $L\in\cL^\ast_F(S_0)$, and by the choice of $G_0$ and $g_0$, we see
that $G_0/L$ is a cyclic $2$-group with $|G_0/L|\ge 4$ and $\<g_0\>+L=G_0$;
furthermore, $S\cap (g_0+L)$ is not contained in a proper coset of $L$, and
$S_0+L=(F\stm G_0)\cup(g_0+L)$, which in view of $S=(G\stm F)\cup S_0$ and
$L\le F$ yields
\begin{equation}\label{e:hren2}
  S+L = (G\stm G_0) \cup (g_0+L).
\end{equation}

It remains to show that $\exp(G/L)=\exp(G_0/L)$ and that $G/G_0$ is
an elementary $2$-group. To prove the former, we observe that
\refe{hren1} and \refe{hren2} yield $2\ast(G/L)\lneqq G_0/L$ and
invoke Lemma \refl{directsum}. To establish the latter, simply
observe that  $2\ast(G/L)\lneqq G_0/L$ implies $2\ast G\le
G_0+L=G_0$, whence $2(g+G_0)=G_0$ for every $g\in G$.
\end{proof}

We can now prove Theorem \reft{main2}.
\begin{proof}[Proof of Theorem \reft{main2}]
We first show that there is at most one subgroup $L\in\cL^\ast_G(S)$ with
\begin{equation}\label{e:main2-0}
  |S+L|-|L|\le |S|-1.
\end{equation}
For a contradiction, assume that $L,\,L'\in\cL^\ast_G(S)$ are distinct, $L$
satisfies \refe{main2-0}, and $|S+L'|-|L'|\le |S|-1$. Find $G_0\le G$ and
$g_0\in G_0$ such that (L1)--(L4) hold, and let $S_0=S\cap(g_0+L)$. It
follows from Lemma \refl{LcapH} that $L'\le G_0$, whence
\begin{equation}\label{e:main2-1}
  |L'|-1 \ge |S+L'|-|S| \ge |S_0+L'|-|S_0|.
\end{equation}
Suppose that $L\nleq L'$ and $L'\nleq L$, and write $t=\pH[L'](S_0)$; that
is, $t$ is the number of $L'$-cosets that intersect $S_0$. Since $S_0$ is not
contained in a proper coset of $L$, and since $L\nleq L'$, we have $t\ge 2$.
Consequently, from $L'\nleq L$ it follows that
  $$ |S_0+L'|-|S_0| \ge t(|L'|-|L\cap L'|) \ge t|L'|/2 \ge |L'|, $$
contradicting \refe{main2-1}. So we may assume either $L\le L'$ or $L'\le
L$; switching the notation, if necessary, and recalling that $L'\ne L$,
we assume that $L<L'$.

Since $L'\in\cL^\ast_G(S)$, there exists a subgroup $G_0'\le G$, lying above
$L'$, and an element $g_0'\in G_0'$ such that $|G_0'|\ge 4|L'|$,
$(S+L')\stm(g_0'+L')=G\stm G_0'$, and $(g_0'+L')\cap S$ is not contained in a
proper coset of $L'$. If $\pH[L'](g_0')=\pH[L'](g_0)$, then $(g_0'+L')\cap
S=(g_0+L')\cap S$, while, in view of $L'\le G_0$, the right-hand side is
contained in an $L$-coset, which, in view of $L<L'$, contradicts that
$(g_0'+L')\cap S$ is not contained in a proper coset of $L'$. Therefore, we
conclude instead that $\pH[L'](g_0)\ne\pH[L'](g_0')$. Thus, since
$|\pi(\pH[L'](S)\stm\{\pH[L'](g_0')\})|=|\pi(G_0'/L')|\ge 4$, it follows from
Proposition \refl{david} that $|\pi(\pH[L'](S)\stm\{\pH[L'](g_0)\})|=1$,
which is equivalent to
  $$ \pi((S+L')\stm(g_0+L'))=L'. $$
Hence, since $L<L'\leq G_0$, so that $(S+L')\stm(g_0+L')=G\stm G_0$, it
follows that $L'=G_0$, whence $S+L'=S+G_0=G$, contradicting the assumption
$L'\in\cL^\ast_G(S)$. This establishes uniqueness of $L\in\cL^\ast_G(S)$
satisfying \refe{main2-0}.

Clearly, Lemma \refl{main} implies assertion (ii) of Theorem 2, and therefore
it remains to establish assertion (i). To this end, suppose that
$L\in\cL^\ast_G(S)$ satisfies \refe{main2-0}, and that $G_0$ and $g_0$ are as
in (L1)--(L4). We will show that $\eta_G(S)\ge|S|$ and
$\kap(\CP(S))=\lam_G(S)=\lam^*_G(S)=|S+L|-|L|$.

Suppose that there exists $H\in\cH_G(S)$ with
\begin{equation}\label{e:main2-2}
  |S+H|-|H| \le |S|-1.
\end{equation}
Then $H\le G_0$ by Lemma \refl{LcapH}. If $H\le L$, then from $(S+2\ast
G)\cap H\neq\est$ we obtain $(S+2\ast G)\cap L\neq\est$, contradicting
(L1)--(L4). Therefore $H\nleq L$.

Let $S_0=(g_0+L)\cap S$, and denote by $t$ the number of $H$-cosets
intersecting $S_0$. In view of \refe{main2-2}, and taking into account
 $H\le G_0$ and $H\nleq L$, we obtain
  $$ |H|-1 \ge |S+H|-|S| \ge |S_0+H|-|S_0|
                                       \ge t(|H|-|H\cap L|) \ge t|H|/2. $$
Hence $t=1$. Thus, since $S_0$ is not contained in a coset of a proper
subgroup of $L$, we conclude that $L\le H$. Consequently, from (L1)--(L3) we
get $2\ast(G/H)=2\ast(G_0/H)$, and thus, in view of
 $(S+2\ast G)\cap H\neq\est$ and taking into account (L4), we have
\begin{equation}\label{e:main2-3}
  \est \ne \pH(S) \cap 2\ast(G/H) = \pH(S) \cap 2\ast(G_0/H) \\
       = \{\pH(g_0)\} \cap 2\ast (G_0/H).
\end{equation}
Since $\pH[L](g_0)$ generates $G_0/L$, it follows from $H\ge L$ that
$\pH(g_0)$ generates the cyclic $2$-group $G_0/H$. Thus \refe{main2-3}
implies that $H=G_0$, whence $S+H=S+G_0=G$, a contradiction. So we conclude
that there are no subgroups $H\in\cH_G(S)$ satisfying \refe{main2-2}; that
is, $\eta_G(S)\ge |S|$. Thus it follows by Lemma \refl{main} that
\begin{equation}\label{e:main2-4}
  \kap(\CP(S))=\min\{\lam^\ast_G(S),|S|\}.
\end{equation}

The uniqueness of $L$, established above, implies that
$\lam^\ast_G(S)=|S+L|-|L|$, and now \refe{main2-0} shows that
  $$ \kap(\CP(S)) \le \lam_G(S) \le \lam^\ast_G(S) = |S+L|-|L| \le |S|-1. $$
Comparing this with \refe{main2-4}, we see that, indeed, the first two
inequalities are actually equalities.
\end{proof}

Finally, we prove Theorem \reft{simple}.
\begin{proof}[Proof of Theorem \reft{simple}]
By Theorem \reft{main2}, we have $\kap(\CP(S))=|S+L|-|L|$ with a subgroup
$L\le G$, belonging to either $\cH_G(S)$ or $\cL^\ast_G(S)$. Let $F\le G$ be
a subgroup that minimizes $|S+F|-|F|$ over all subgroups with $S+F\ne G$.
Assuming that
\begin{equation}\label{e:simple1}
  |S+F|-|F| < |S+L|-|L| \le |S|-1,
\end{equation}
we will obtain a contradiction; evidently, this will prove the assertion.

>From Lemma \refl{geom} and \refe{simple1}, it follows that either $S+F+L=S+L$
or $S+F+L=S+F$; in either case,
\begin{equation}\label{e:simple2}
  S+F+L \ne G.
\end{equation}

Suppose first that $|L|\le |F|$. Then Lemma \refl{geom} yields $S+F+L=S+F$,
and thus
  $$ |S+F+L|-|F+L| = |S+F|-|F+L|. $$
The minimality of $F$ now implies that $|F+L|=|F|$, whence $L\le F$. If
$L\in\cH_G(S)$, then it follows in view of $L\le F$ and $S+F\ne G$ that
$F\in\cH_G(S)$, implying $\kap(\CP(S))\le|S+F|-|F|$. However, since
$\kappa(\Cay_G^+(S))=|S+L|-|L|$, this contradicts \refe{simple1}. Therefore
we may assume $L\in\cL^\ast_G(S)$. Let $G_0$ be the subgroup from the
definition of $\cL^\ast_G(S)$. By Lemma \refl{LcapH} we then have $L\le F\le
G_0$, whence
  $$ |S+F| = |G\stm G_0| + |F| = (|S+L|-|L|) + |F|, $$
which contradicts \refe{simple1} once more.

Next, suppose that $|F|\le|L|$. Thus it follows by Lemma \refl{geom} that
$S+L=S+F+L$. Hence
\begin{equation}\label{e:simple3}
  |S+F+L|-|F+L| = |S+L|-|F+L|.
\end{equation}

If $L\in\cH_G(S)$, then it follows in view of $L\le F+L$ and \refe{simple2}
that $F+L\in \cH_G(S)$; now \refe{simple3} and the minimality of $L$ give
$|F+L|=|L|$, leading to $F\le L$. We proceed to show that this holds in the
case $L\in\cL^\ast_G(S)$ as well. In this case, in view of \refe{simple3} and
\refe{simple1}, Lemma \refl{LcapH} gives $F+L\le G_0$, where $G_0$ is the
subgroup from the definition of $\cL^\ast_G(S)$. Thus (as in the previous
paragraph)
  $$ |S+F+L| = |G\stm G_0| + |F+L| = (|S+L|-|L|) + |F+L|. $$
Hence, since $|S+F+L|=|S+L|$, we obtain $|F+L|=|L|$, and therefore $F\le L$,
as desired.

We have just shown that $F\le L$ holds true in either case. Consequently,
from $|S+L|-|L|<|S|\le |S+F|$ and divisibility considerations, it follows
that indeed $|S+L|-|L|\le|S+F|-|F|$, contradicting \refe{simple1} and
completing the proof.
\end{proof}

\bigskip


\begin{thebibliography}{CGW03}

\bibitem[A]{b:a}
  {\sc N.~Alon},
  Large sets in finite fields are sumsets,
  \emph{J.~Number Theory}, to appear.

\bibitem[CGW03]{b:cgw}
  {\sc B.~Cheyne, V.~Gupta}, and {\sc C.~Wheeler},
  Hamilton Cycles in Addition Graphs,
  \emph{Rose-Hulman Undergraduate Math. Journal} {\bf 1} (4) (2003)
  (electronic).


\bibitem[C92]{b:chu}
  {\sc F.R.K.~Chung},
  Diameters and eigenvalues,
  \emph{J.~Amer. Math. Soc.} {\bf 2} (2) (1989), 187--196.

\bibitem[Gr05]{b:g}
  {\sc B.J.~Green},
  Counting sets with small sumset, and the clique number of random Cayley
    graphs,
  \emph{Combinatorica} {\bf 25} (2005), 307--326.

\bibitem[Gk01]{b:d1}
  {\sc D.J.~Grynkiewicz},
  Quasi-periodic decompositions and the Kemperman structure theorem,
  \emph{European Journal of Combinatorics} {\bf 26} (5) (2005), 559--575.

\bibitem[Gk02]{b:d2}
  {\sc \bysame},
  Sumsets, Zero-Sums and Extremal Combinatorics,
  \emph{Ph.D. Dissertation}, Caltech (2005).

\bibitem[Km60]{b:k}
  {\sc J.H.B.~Kemperman},
  On small subsets in an abelian group,
  \emph{Acta Mathematica} {\bf 103} (1960), 63--88.

\bibitem[Kn53]{b:kn1}
  {\sc M.~Kneser},
  Absch\"atzung der asymptotischen Dichte von Summenmengen,
  \emph{Math. Z.} {\bf 58} (1953), 459--484.

\bibitem[Kn55]{b:kn2}
  \bysame,
  Ein Satz \"uber abelsche Gruppen mit Anwendungen auf die
  Geometrie der Zahlen,
  \emph{Math. Z.} {\bf 61} (1955), 429--434.

\bibitem[L05]{b:l1}
  {\sc V.F.~Lev},
  Critical pairs in abelian groups and Kemperman's structure theorem,
  \emph{Internation Journal of Number Theory} {\bf 2} (3) (2006), 379--396.

\bibitem[L]{b:l2}
  {\sc \bysame},
  Sums and differences along Hamiltonian cycles,
  \emph{Submitted}.


\bibitem[Mn76]{b:m}
  {\sc H.B.~Mann},
  Addition theorems: the addition theorems of group theory and number theory.
  Reprint, with corrections, of the 1965 original.
  Robert E. Krieger Publishing Co., Huntington, N.Y., 1976.

\bibitem[Mr83]{b:mr}
  {\sc D.~Maru\v si\v c},
  Hamiltonian circuits in Cayley graphs,
  \emph{Discrete Math.} {\bf 46} (1) (1983), 49--54.

\end{thebibliography}
\end{document}